# On the isometrization of groups of homeomorphisms

by Fredric D. Ancel

**Abstract.** Let G be a group of homeomorphisms of a topological space X. G is *(properly) isometrizable* if there exists a G-invariant (proper) gauge structure on X. G is *equiregular* if for every x ∈ X and every open neighborhood U of x in X there is an open neighborhood V of x in X such that cl(V) ⊂ U and every y ∈ X has an open neighborhood $N_y$ with the property that for every g ∈ G, if g($N_y$) ∩ cl(V) ≠ ∅, then g($N_y$) ⊂ U. G is *nearly proper* if for all compact subsets A and B ⊂ X, cl(∪ { g(A) : g ∈ G and g(A) ∩ B ≠ ∅ }) is compact.

**The Isometrization Theorem.** If X is a Hausdorff space and G\X is a paracompact regular space, then: G is isometrizable if and only if G is equiregular.

**The Proper Isometrization Theorem.** If X is a locally compact σ-compact Hausdorff space and G\X is a regular space, then: G is properly isometrizable if and only if G is equiregular and nearly proper.

G *acts properly on* X if for all compact subsets A and B of X, the subset $G_{A,B}$ = { g ∈ G: g(A) ∩ B ≠ ∅ } is compact when G is endowed with the compact-open topology. The Proper Isometrization Theorem has the following corollary.

**Theorem of Abel-Manoussos-Noskov [AMN].** If X is a locally compact σ-compact Hausdorff space and G acts properly on X, then X is properly isometrizable.

## 1. Introduction

This article is, to some extent, a re-exposition of the paper [AMN] with two key lemmas added that make it possible to generalize the results of [AMN]. I was led to read this paper by a (still-unanswered) question asked me by Craig Guilbault. His question is recorded near the end of the introduction.

We begin by setting the context with some definitions and comments about them, followed by statements of our main results.

**Notation.** Let X be a topological space and let ρ : X × X → [0,∞] be a function. For x ∈ X and є ∈ (0,∞], let $\mathcal{N}_\rho$(x,є) = { y ∈ X : ρ(x,y) < є } and let $\mathcal{N}_\rho$(x,є] = { y ∈ X : ρ(x,y) ≤ є }. For A ⊂ X and є ∈ (0,∞], let $\mathcal{N}_\rho$(A,є) = ∪$_{x \in A}\mathcal{N}_\rho$(x,є).

We give a definition of a *gauge structure* on a topological space that differs slightly from the standard definition in that it reflects the given topology on the underlying space.

**Definition.** Let X be a topological space. A function ρ : X × X → [0,∞) is a *gauge* or *pseudometric* on X if it satisfies the following four conditions.

*2020 Mathematics Subject Classification:* 22F05, 54H15, 57S05

*1)* $\rho(x,x) = 0$ for all $x \in X$.

*2)* $\rho(x,y) = \rho(y,x)$ for all $x, y \in X$.

*3)* $\rho(x,z) \leq \rho(x,y) + \rho(y,z)$ for all $x, y, z \in X$.

*4)* $\mathcal{N}_\rho(x,\epsilon)$ is an open subset of X for all $x \in X$ and all $\epsilon \in (0,\infty)$.

(The fourth condition is not part of the usual definition of a gauge.)
$\rho$ is a *metric* on X if it also satisfies:

*5)* $\rho(x,y) = 0 \Rightarrow x = y$ for all $x, y \in X$.

Let $\mathcal{P}$ be a collection of gauges on X. Then { $\mathcal{N}_\rho(x,\epsilon) : \rho \in \mathcal{P}$, $x \in X$, $\epsilon \in (0,\infty)$ } is a subbasis for a topology on X called the *topology determined by $\mathcal{P}$*. (This topology is clearly a subset of the given topology on X.) Two collections of gauges on X are *equivalent* if they determine the same topology on X. $\mathcal{P}$ is *separating* if for all distinct x, $y \in X$, there is a $\rho \in \mathcal{P}$ such that $\rho(x,y) > 0$. (Thus, $\mathcal{P}$ is separating if and only if the topology determined by $\mathcal{P}$ is Hausdorff.) If $\mathcal{P}$ is separating and $\mathcal{P}$ determines the given topology on X, then $\mathcal{P}$ is called a *gauge structure on X* . A topological space that has a gauge structure is called a *gauge space.*

We refer to pages 198 to 200 of [D] as a source of information about gauge spaces. To provide perspective, we remind the reader of a theorem found in these pages of [D] that characterizes gauge spaces: for a Hausdorff space X, the following conditions are equivalent:
- X is a gauge space.
- X is completely regular.
- There is a set J such that X embeds in $[0,1]^J$.

**Definition.** Let G be a group of homeomorphisms of a topological space X. Let $\mathcal{P}$ be a collection of gauges on X. If for all x, $y \in X$, all $g \in G$ and all $\rho \in \mathcal{P}$, $\rho(g(x),g(y)) = \rho(x,y)$, then we call G a *$\mathcal{P}$-isometry group* and we say that $\mathcal{P}$ is *G-invariant*. G is *isometrizable* if X has a G-invariant gauge structure. G is *singly isometrizable* if X has a one-element G-invariant gauge structure {$\rho$}. (In other words, G is singly isometrizable if and only if there is a G-invariant metric $\rho$ on X that determines the given topology on X.)

**Definition.** A group G of homeomorphisms of a topological space X is *equiregular* if for every $x \in X$ and every open neighborhood U of x in X there is an open neighborhood V of x in X such that cl(V) $\subset$ U and every $y \in X$ has an open neighborhood $N_y$ with the property that for every $g \in G$, if $g(N_y) \cap cl(V) \neq \emptyset$, then $g(N_y) \subset U$.

**Remark.** Observe that if G is an equiregular group of homeomorphisms of a space X, then X is a regular space. In addition, the equiregularity of G entails that the action of G on X "witnesses" the regularity of X.

We now state one of the two main results of this article.

**Note.** In this article we *don't* assume that paracompact spaces are *Hausdorff*, only that in a paracompact space, every open cover has a locally finite refinement.

**The Isometrization Theorem 1.1.** If X is a Hausdorff space, G is a group of homeomorphisms of X, and G\X is a paracompact regular space, then: G is isometrizable if and only if G is equiregular.

**Corollary to the Isometrization Theorem 1.2.** Suppose X is a separable metrizable space, G is a group of homeomorphisms of X, and G\X is a paracompact regular space. If G is equiregular, then G is singly isometrizable.

We sketch the proof of this corollary after we introduce a bit of terminology.

**Definition.** If $\rho$ is a gauge on a topological space X, define the *decapitation* $\overline{\rho}$ of $\rho$ by

$$\overline{\rho}(x,y) = \min \{ \rho(x,y),\ 1 \} \text{ for x, y} \in X.$$

Then $\overline{\rho}$ is a gauge on X that is equivalent to $\rho$.

**Definition.** If $\rho_1, \rho_2, \cdots, \rho_k$ are gauges on a topological space X, define *max{$\rho_1,\rho_2,\cdots,\rho_k$}* by

$$max\{\rho_1,\rho_2,\cdots,\rho_k\}(x,y) = \max \{ \rho_1(x,y),\ \rho_2(x,y),\ \cdots,\ \rho_k(x,y) \} \text{ for x, y} \in X.$$

Then max{$\rho_1,\rho_2,\cdots,\rho_k$} is a gauge on X that is equivalent to { $\rho_1, \rho_2, \cdots, \rho_k$ }. Moreover, for $x \in X$ and $\epsilon \in (0,\infty)$, $\mathcal{N}_{max\{\rho_1,\rho_2,\cdots,\rho_k\}}(x,\epsilon) = \bigcap_{1 \le i \le k} \mathcal{N}_{\rho_i}(x,\epsilon)$.

**Definition.** Let $\mathcal{P}$ be a collection of gauges on a topological space X. $\mathcal{P}$ is *closed under maximization* of $\rho_1, \rho_2, \cdots, \rho_k \in \mathcal{P} \Rightarrow max\{\rho_1,\rho_2,\cdots,\rho_k\} \in \mathcal{P}$. Let

$$max(\mathcal{P}) = \{ max\{\rho_1,\rho_2,\cdots,\rho_k\} : \rho_1, \rho_2, \cdots, \rho_k \in \mathcal{P} \}.$$

**Observations.** Let $\mathcal{P}$ be a collection of gauges on a topological space X.
- If $\mathcal{P}$ is a gauge structure on X, then so is max($\mathcal{P}$).
- max($\mathcal{P}$) is closed under maximization.
- If $\mathcal{P}$ is a gauge structure on X that is closed under maximization, then { $\mathcal{N}_\rho(x,\epsilon) : \rho \in \mathcal{P}, x \in X, \epsilon \in (0,\infty)$ } is a basis (not merely a subbasis) for the topology on X.

**Further observations.** Let G be a group of homeomorphisms of a topological space X.
- If $\rho$ is a G-invariant gauge on X, then $\overline{\rho}$ is G-invariant.
- If $\rho_1, \rho_2, \cdots, \rho_k$ are G-invariant gauges on X, then max{$\rho_1,\rho_2,\cdots,\rho_k$} is G-invariant.
- Hence, if $\mathcal{P}$ is a G-invariant gauge structure on X, then max($\mathcal{P}$) is G-invariant.

**Proof of the Corollary to the Isometrization Theorem.** The Isometrization Theorem provides a G-invariant gauge structure $\mathcal{P}$ on X which we can assume is closed under maximization. Since X has a countable basis, then every basis for X contains a countable basis. Therefore, X has a basis of the form { $\mathcal{N}_{\rho_i}(x_i, \epsilon_i)$: i ≥ 1 } where $\rho_i \in \mathcal{P}$, $x_i \in X$ and $\epsilon_i \in (0, \infty)$ for each i ≥ 1. Hence, { $\rho_i$ : i ≥ 1 } is a countable G-invariant gauge structure on X. Consequently, { $\bar{\rho}_i$ : i ≥ 1 } is a countable G-invariant gauge structure on X. Define σ : X × X → [0,∞) by σ(x,y) = sup { $2^{-i}\bar{\rho}_i(x,y)$ : i ≥ 1 }. Then σ is a G-invariant metric on X that is equivalent to { $\bar{\rho}_i$ : i ≥ 1 }. Hence, σ determines the given topology on X. This proves G is singly isometrizable. ∎

**Question.** Suppose X is a non-separable metrizable space, G is a group of homeomorphisms of X, and G\X is a paracompact regular space. If G is equiregular, must G be singly isometrizable?

**Definition.** Let X be a topological space. A gauge ρ on X is *proper* if cl($\mathcal{N}_\rho(x, \epsilon)$) is compact whenever x ∈ X and ∈ ∈ (0,∞). A gauge structure $\mathcal{P}$ on X is *proper* if every element of $\mathcal{P}$ is proper. A topological space that has a proper gauge structure is called a *proper gauge space.*

**Observation.** If a topological space X has a proper gauge ρ, then X is locally compact and σ-compact. σ-compactness follows from the fact that X = $\bigcup_{n \geq 1}$cl($\mathcal{N}_\rho(x,n)$).

**Definition.** Let G be a group of homeomorphisms of a topological space X. G is *properly isometrizable* if X has a G-invariant proper gauge structure. G is *singly properly isometrizable* if X has a one-element G-invariant proper gauge structure {ρ} (i.e., ρ is a G-invariant proper metric on X that determines the given topology on X.)

**Definition.** A group G of homeomorphisms of a space X is *nearly proper* if for all compact subsets A and B of X, cl($\bigcup$ { g(A) : g ∈ G and g(A) ∩ B ≠ ∅ }) is compact.

We now state the second main results of this article.

**The Proper Isometrization Theorem 1.3.** If X is a locally compact σ-compact Hausdorff space and G\X is a regular space, then: G is properly isometrizable if and only if G is equiregular and nearly proper.

**Corollary to the Proper Isometrization Theorem 1.4.** Suppose X is a locally compact σ-compact metrizable space and G\X is a regular space. If G is equiregular and nearly compact, then G is singly properly isometrizable.

Before proving this corollary, we make two observations.

**Observations.** Let X be a topological space.

- If ρ and σ are gauges on X and σ is proper, then max{ρ,σ} is a proper gauge on X that is equivalent to { ρ, σ }.
- Hence, if $\mathcal{P}$ is a gauge structure on X and σ is a proper gauge on X, then { max{ρ,σ} : ρ ∈ $\mathcal{P}$ } is a proper gauge structure on X.

**Proof of the Corollary to the Proper Isometrization Theorem.** Since X is a locally compact σ-compact metrizable space, then X is separable. Also the orbit map π : X → G\X, being an open map, preserves local compactness and σ-compactness. It follows that G\X is paracompact. (The proofs of Proposition 5.2 and Lemma 5.3 give a detailed justification of these assertions.) Therefore, the Corollary to the Isometrization Theorem provides a G-invariant metric ρ on X that determines the given topology on X. Also, the Proper Isometrization Theorem tells us that there is a G-invariant proper gauge σ on X. Hence, max{ρ,σ} is a G-invariant proper metric on X that determines the given topology. ∎

Next we present two examples which show that equiregularity and near properness are independent notions.

**Example 1.1.** Let X = { 0, 1, 2, ··· } with the discrete metric ρ. (ρ(x,y) = 1 if x ≠ y.) Let G be the homeomorphism group of X. Since ρ is G-invariant, then G is isometrizable and, hence, equiregular by the Isometrization Theorem. However, G is not nearly proper. Indeed, if C = { 0, 1 } and for each n ≥ 1, $g_n$ : X → X is the homeomorphism that transposes 1 and n and fixes the other points of X, then $g_n(C) \cap C \neq \emptyset$ for each n ≥ 1 and $\bigcup_{n \geq 1} g_n(C) = X$. Hence, G is singly isometrizable but not properly isometrizable.

**Example 1.2.** Let X = [0,1]. Define the homeomorphism g : X → X by $g(x) = x^2$. Let G be the group of homeomorphisms of X generated by g. Since X is compact, G is nearly compact. However, for every open neighborhood N of 1 in X, $\bigcup_{n \geq 1} g^{(n)}(N) = (0,1]$. Hence, G is not equiregular. Indeed, for every open neighborhood N of 1 in X, there is an n ≥ 1 such that $1 \in g^{(n)}(N)$ and $g^{(n)}(N) \not\subset (½,1]$.

**Definition.** A group G of homeomorphisms of a topological space X *acts properly on* X if for all compact subsets A and B of X, the subset $G_{A,B}$ = { g ∈ G: g(A) ∩ B ≠ ∅ } is compact when G is endowed with the compact-open topology.

There are several isometrization theorems for properly acting groups of homeomorphisms in the published literature including [AdN] and [AMN]. The strongest of these theorems are:

**Theorem 1.5. [AMN].** If G is a group of homeomorphisms of a locally compact σ-compact Hausdorff space X and G acts properly on X, then G is properly isometrizable.

**Corollary 1.6. [AMN].** If G is a group of homeomorphisms of a locally compact σ-compact metrizable space X and G acts properly on X, then G is singly properly isometrizable.

These results are corollaries to the previously stated Proper Isometrization Theorem and its corollary. Indeed, if G is a properly acting group of homeomorphisms of a locally compact σ-compact Hausdorff space X, then relatively simple arguments reveal that: *1)* G\X is locally compact and Hausdorff, hence, regular; *2)* G is equiregular; and *3)* G is nearly proper. Therefore, the Proper Isometrization Theorem and its corollary yield the theorem and corollary of [AMN]. The following example shows that the Proper Isometrization Theorem and its corollary apply to a wider class of homeomorphism groups than do the corresponding results of [AMN].

**Example 1.3.** For every a ∈ ℚ = { rational numbers }, define the homeomorphism $g_a$ : ℝ → ℝ by $f_a(x) = x + a$. Let G = { $g_a$ : a ∈ ℚ }. Then G is a group of homeomorphisms of ℝ and the usual metric on ℝ is a proper G-invariant metric. Thus, the Proper Isometrization Theorem applies to G and tells us that G is equiregular and nearly proper. However, since G does not act properly on ℝ, then G is outside the range of application of the results of [AMN]. To see that the action of G on ℝ is not proper, first note that if G is endowed with the compact-open topology, then a ↦ $g_a$ : ℚ → G is a homeomorphism. Hence, the set $G_{\{0\},[0,1]}$ = { $g_a$ : a ∈ ℚ ∩ [0,1] } is non-compact.

Next we state a currently unanswered question asked by Craig Guilbault that motivated the author to study the results in [AdM] and [AMN].

**Definition.** Let ρ be a metric on a topological space X. A map γ : [a,b] → X is a *ρ-geodesic* if ρ(γ(s),γ(t)) = l s − t l for all s, t ∈ [a,b]. ρ is a *geodesic metric* if for all x, y ∈ X, there is a ρ-geodesic γ : [a,b] → X such that γ(a) = x and γ(b) = y.

The Bing-Moise Convex Metrization Theorem [Bi], [Mo] implies that every locally compact, locally connected, connected metrizable space has a proper geodesic metric.

**Question (Craig Guilbault).** If G is an equiregular and nearly proper group of homeomorphisms of a locally compact, locally connected, connected metrizable space X, does X have a G-invariant proper ***geodesic*** metric?

**Outline of remaining article.** In section 2, we describe the Alexandroff-Urysohn methods for constructing general gauge structures on a topological space. In section 3, focus on the construction of proper gauge structures. In section 4, we show to make the gauge structures and proper gauge structures of sections 2 and 3 invariant with respect to a group of homeomorphisms of the underlying space. In section 5, we state and prove two key lemmas. These lemmas allow us to prove the Isometrization Theorem and Proper Isometrization Theorem by showing that the hypotheses of these theorems are sufficient to carry out the construction of invariant gauge structures and proper gauge structures of the sort described in section 3. In section 6, we show how to deduce the theorems of [AMN] from the Proper Isometrization Theorem.

**Acknowledgements.** I acknowledge Craig Guilbault's question for focusing my attention on the article [AMN]. This article was in turn my primary source for background and ideas leading to the results presented here. I also credit the article [R] written when both its author and I were graduate students in Madison Wisconsin for making me aware of the Alexandroff-Urysohn metrization method.

## 2. The construction of gauge structures

There are two steps to the process we describe here. The first step produces gauges that may be infinite valued. The second step modifies gauges to remove infinite values. We begin by giving the formal definition of an infinite-valued gauge.

**Definition.** Let X be a topological space. A function $\rho : X \times X \to [0,\infty]$ is an *infinite-valued gauge* or, more briefly, an *$\infty$-gauge* on X if it satisfies conditions *1), 2), 3)* and *4)* stated in the definition of "gauge". If $\rho$ also satisfies condition *5)*, it is an *$\infty$-metric* on X. Clearly, every gauge on X is an $\infty$-gauge, and every metric is an $\infty$-metric. Let $\mathcal{P}$ be a collection of $\infty$-gauges on X. Then { $\mathcal{N}_\rho(x,\epsilon) : \rho \in \mathcal{P}, x \in X, \epsilon \in (0,\infty)$ } is a subbasis for a topology on X called the *topology determined by $\mathcal{P}$*. (This topology is a subset of the given topology on X.) Two collections of $\infty$-gauges on X are *equivalent* if they determine the same topology on X. $\mathcal{P}$ is *separating* if for all distinct x, y $\in$ X, there is a $\rho \in \mathcal{P}$ such that $\rho(x,y) > 0$. (Thus, $\mathcal{P}$ is separating if and only if the topology determined by $\mathcal{P}$ is Hausdorff.) If $\mathcal{P}$ is separating and $\mathcal{P}$ determines the given topology on X, then $\mathcal{P}$ is called an *$\infty$-gauge structure on X*. A topological space that has an $\infty$-gauge structure is called an *$\infty$-gauge space.*

**Definition.** Suppose X is a topological space and $\rho$ is an $\infty$-gauge on X. For x $\in$ X, the set $\mathcal{N}_\rho(x,\infty) = \{ y \in X : \rho(x,y) < \infty \}$ is called the *$\rho$-crevasse of x.*

Observe that the collection of $\rho$-crevasses of the points of a space X forms a partition of X into non-empty open and, hence, closed subsets.

**Definition.** If $\rho$ is an $\infty$-gauge on a topological space X, then, as before, define the *decapitation $\overline{\rho}$* of $\rho$ by

$$\overline{\rho}(x,y) = \min \{ \rho(x,y), 1 \} \text{ for x, y} \in X.$$

Then $\overline{\rho}$ is a (finite-valued) gauge on X that is equivalent to $\rho$.

**Definition.** Let $\mathcal{U}$ be an open cover of a topological space X and let A $\subset$ X. Let

$$Star(A, \mathcal{U}) = \cup \{ U \in \mathcal{U} : A \cap U \neq \emptyset \}.$$

For n ≥ 1, inductively define *Star$^n$(A,$\mathcal{U}$)* by

$$Star^1(A,\mathcal{U}) = \text{Star}(A,\mathcal{U}) \text{ and } Star^{n+1}(A,\mathcal{U}) = \text{Star}(\text{Star}^n(A,\mathcal{U}),\mathcal{U}).$$

**Definition.** Let $\mathcal{U}$ and $\mathcal{V}$ be open covers of a topological space X. We say that $\mathcal{V}$ *star-refines* $\mathcal{U}$ if { Star(x,$\mathcal{V}$) : x ∈ X } refines $\mathcal{U}$.

**Example 2.1.** If ρ is an ∞-gauge on a topological space X, then for every ϵ ∈ (0,∞), { $\mathcal{N}_\rho$(x,ϵ/2) : x ∈ X } star-refines { $\mathcal{N}_\rho$(x,ϵ) : x ∈ X }.

**Definition.** Let $\mathcal{U}$ be an open cover of a topological space X. A k-tuple $(U_1,U_2,\cdots,U_k)$ of elements of $\mathcal{U}$ such that $U_i \cap U_{i+1} \neq \emptyset$ for $1 \leq i < k$ is called a *$\mathcal{U}$-chain with k links*. If $(U_1,U_2,\cdots,U_k)$ is a $\mathcal{U}$-chain such that x ∈ $U_1$ and y ∈ $U_k$, then we say that $(U_1,U_2,\cdots,U_k)$ *joins x to y.* For x, y ∈ X, let $\mathcal{U}$(x,y) denote the set of all $\mathcal{U}$-chains joining x to y.

**Definition.** Let $\mathcal{U}$ be an open cover of a topological space X. For each x ∈ X, call { y ∈ X : $\mathcal{U}$(x,y) ≠ ∅ } the *$\mathcal{U}$-component of x.* If for all x, y ∈ X, $\mathcal{U}$(x,y) ≠ ∅, we say X is *$\mathcal{U}$-connected.*

Let $\mathcal{U}$ be an open cover of a topological space X. Observe that the collection of $\mathcal{U}$-components forms a partition of X into non-empty open and, hence, closed subsets. Also, observe that X is connected if and only if X is $\mathcal{U}$-connected for every open cover $\mathcal{U}$ of X.

**Definition.** Let $\mathcal{U}$ and $\mathcal{V}$ be open covers of a topological space X and let k ≥ 1. We say that $\mathcal{V}$ *k-refines* $\mathcal{U}$ if every $\mathcal{V}$-chain with k or fewer links is contained in an element of $\mathcal{U}$.

Let $\mathcal{U}$, $\mathcal{V}$ and $\mathcal{W}$ be open covers of a topological space X. Observe that if $\mathcal{V}$ star-refines $\mathcal{U}$, then $\mathcal{V}$ 2-refines $\mathcal{U}$. Further observe that if $\mathcal{V}$ 2-refines $\mathcal{U}$ and $\mathcal{W}$ 2-refines $\mathcal{V}$, then $\mathcal{W}$ 4-refines $\mathcal{U}$; and if $\mathcal{W}$ 4-refines $\mathcal{U}$, then $\mathcal{W}$ 3-refines $\mathcal{U}$.

**Definition.** If { $\mathcal{U}_n$ : n ≥ 1 } is a sequence of open covers of a topological space X such that { Star(x,$\mathcal{U}_n$) : x ∈ X and n ≥ 1 } is a basis for a topology on X, then this topology is called the *topology determined by* { $\mathcal{U}_n$ : n ≥ 1 }, and { $\mathcal{U}_n$ : n ≥ 1 } is called a *development in* X.

We observe that the requirement that { Star(x,$\mathcal{U}_n$) : x ∈ X and n ≥ 1 } be a basis for a topology on the space X is equivalent to requiring that for all x, y ∈ X and all m, n ≥ 1, if z ∈ Star(x,$\mathcal{U}_m$) ∩ Star(y,$\mathcal{U}_n$), then there is an r ≥ 1 such that Star(z,$\mathcal{U}_r$) ⊂ Star(x,$\mathcal{U}_m$) ∩ Star(y,$\mathcal{U}_n$). We remark the usual definition of a development requires { Star(x,$\mathcal{U}_n$) : x ∈ X and n ≥ 1 } to be a basis for the given topology on X. We have weakened this requirement so that the topology determined by { $\mathcal{U}_n$ : n ≥ 1 } is allowed to be a *subset* of the given topology on X.

**Definition.** Let $\{\,\mathcal{U}_n : n \geq 1\,\}$ be a development in a topological space X. $\{\,\mathcal{U}_n : n \geq 1\,\}$ is a *star-development* if $\forall\, n \geq 1$, $\mathcal{U}_{n+1}$ star-refines $\mathcal{U}_n$. For $k \geq 1$, $\{\,\mathcal{U}_n : n \geq 1\,\}$ is a *k-development* if $\forall\, n \geq 1$, $\mathcal{U}_{n+1}$ k-refines $\mathcal{U}_n$.

Let $\{\,\mathcal{U}_n : n \geq 1\,\}$ be a development in a topological space X. Observe that if $\{\,\mathcal{U}_n : n \geq 1\,\}$ is a star-development then $\{\,\mathcal{U}_n : n \geq 1\,\}$ is a 2-development. Further observe that if $\{\,\mathcal{U}_n : n \geq 1\,\}$ is a 2-development, then $\{\,\mathcal{U}_{2n-1} : n \geq 1\,\}$ is a 4-development and, hence, a 3-development; and in this situation, $\{\,\mathcal{U}_n : n \geq 1\,\}$ and $\{\,\mathcal{U}_{2n-1} : n \geq 1\,\}$ determine the same topology on X.

**Example 2.2.** Let $\rho$ be an $\infty$-gauge on a topological space X. For every $n \geq 1$, let $\mathcal{U}_n = \{\,\mathcal{N}_\rho(x, 2^{-n}) : x \in X\,\}$. Then $\{\,\mathcal{U}_n : n \geq 1\,\}$ is a star-development in X. Observe that for every $x \in X$ and every $n \geq 1$, $\mathcal{N}_\rho(x, 2^{-n-1}) \subset \text{Star}(x, \mathcal{U}_{n+1}) \subset \mathcal{N}_\rho(x, 2^{-n})$. Hence, the topologies on X determined by $\rho$ and by $\{\,\mathcal{U}_n : n \geq 1\,\}$ coincide.

**Definition.** If $\mathbb{D}$ is a collection of developments in a topological space X, then $\{\,\text{Star}(x, \mathcal{U}_n) : \{\,\mathcal{U}_n : n \geq 1\,\} \in \mathbb{D},\ x \in X \text{ and } n \geq 1\,\}$ is a subbasis for a topology on X called the *topology determined by* $\mathbb{D}$.

**The Alexandroff-Urysohn $\infty$-gauge Construction Theorem 2.1.** If $\{\,\mathcal{U}_n : n \geq 1\,\}$ is a 3-development in a topological space X, then there is an $\infty$-gauge $\rho$ on X, called the *A-U distance function associated with* $\{\,\mathcal{U}_n : n \geq 1\,\}$, with the property that for every $x \in X$ and every $n \geq 1$, $\mathcal{N}_\rho(x, 2^{-n}) \subset \text{Star}(x, \mathcal{U}_n) \subset \mathcal{N}_\rho(x, 2^{-n})$. Hence, $\rho$ and $\{\,\mathcal{U}_n : n \geq 1\,\}$ determine the same topology on X. $\rho$ is defined from $\{\,\mathcal{U}_n : n \geq 1\,\}$ according to the following recipe.

Let $\mathcal{U}^* = \cup_{n \geq 1} \mathcal{U}_n$. Define $\delta : \mathcal{U}^* \to \{0\} \cup \{\,2^{-n} : n \geq 1\,\}$ by

$$\delta(U) \ = \ \inf\{\,2^{-n} : U \in \mathcal{U}_n\,\}.$$

If $C = (U_1, U_2, \cdots, U_k)$ is a $\mathcal{U}^*$-chain, define $\lambda(C) \in [0, \infty)$ by

$$\lambda(C) \ = \ \Sigma_{1 \leq i \leq k}\, \delta(U_i).$$

Define the *A-U distance function* $\rho : X \times X \to [0, \infty]$ *associated with* $\{\,\mathcal{U}_n : n \geq 1\,\}$ by

$$\rho(x,y) \ = \ \inf\left(\{\,\lambda(C) : C \in \mathcal{U}^*(x,y)\,\} \cup \{\infty\}\right).$$

Thus, the set of $\rho$-crevasses equals the set of $\mathcal{U}_1$-crevasses.

Observe that $\rho(x,y)$ is the infimum of the lengths of $\mathcal{U}^*$-chains joining x to y.

A proof of Theorem 2.1 can be obtained via minor changes to the proof of the classical Alexandroff-Urysohn metrization theorem found in the literature. For this, we refer the reader to either the original paper of Alexandroff and Urysohn [AU] or to Rolfsen's article [R]. In [AU], it is hypothesized that the development $\{\,\mathcal{U}_n : n \geq 1\,\}$ determines the given topology on X which was assumed to be Hausdorff. It was also assumed that $\mathcal{U}_1 = \{X\}$. Therefore, the original theorem yields a finite-valued metric on X which determines the given topology. In Theorem 2.1, weaker hypotheses are assumed resulting in weaker conclusions: ρ is merely an ∞-gauge on X which determines a topology on X that is a subset of the given topology. (The Alexandroff-Urysohn proof was appropriated by A. Weil for his proof that uniform spaces have gauge structures. There are many expositions of Weil's proof in the topology literature. For instance see Theorem 11.4 on pages 203-204 of [D]. Each of these proofs can easily be modified to provide a proof of Theorem 2.1.)

For the record, we state the original 1923 theorem of Alexandroff and Urysohn.

**The Classical Alexandroff-Urysohn Metrization Theorem [AU].** If a Hausdorff space X has a 2-development that determines its topology, then X is metrizable.

If ρ is the A-U distance function associated with a 3-development { $\mathcal{U}_n$ : n ≥ 1 } in a topological space X, then the decapitation of ρ is a finite-valued gauge on X that is equivalent to ρ. This observation gives us the following two corollaries.

**Corollary 2.2.** If { $\mathcal{U}_n$ : n ≥ 1 } is a 3-development in a space X, then the decapitation of the associated A-U distance function is a (finite-valued) gauge on X that determines the same topology as { $\mathcal{U}_n$ : n ≥ 1 }. ∎

**Corollary 2.3.** If $\mathbb{D}$ is a collection of 3-developments in a Hausdorff space X that determines the given topology on X, then the collection of all the decapitations of the A-U distance functions associated with the elements of $\mathbb{D}$ is a gauge structure on X that determines the given topology on X. ∎

## 3. The construction of proper gauge structures

**Definition.** Let X be a topological space. An ∞-gauge ρ on X is *proper* if cl($\mathcal{N}_\rho(x,\epsilon)$) is compact whenever x ∈ X and $\epsilon$ ∈ (0,∞). An ∞-gauge structure $\mathcal{P}$ on X is *proper* if every element of $\mathcal{P}$ is proper. A topological space that has a proper ∞-gauge structure is called a *proper ∞-gauge space.*

Observe that if ρ is a proper ∞-gauge on X, then X is locally compact because $\mathcal{N}_\rho(x,1)$ is an open neighborhood of x with compact closure for every x ∈ X. We have already noted that if ρ is a proper finite-valued gauge on X, then X is σ-compact. Hence, if ρ is an ∞-gauge on X that is equivalent to a proper finite-valued gauge, then X is σ-compact.

**Lemma 3.1.** If ρ is a proper ∞-gauge on a topological space X, then for every compact subset A of X and every $\epsilon$ ∈ (0,∞), cl($\mathcal{N}_\rho(A,\epsilon)$) is compact.

**Proof.** Since the ρ-crevasses form an open cover of X, there is a finite sequence $Y_1$, $Y_2$, ⋯ , $Y_k$ of ρ-crevasses such that A ⊂ $\cup_{1 \leq i \leq k} Y_i$. For 1 ≤ i ≤ k, let $A_i$ = A ∩ $Y_i$. Since each $Y_i$ is a closed subset of X, then each $A_i$ is compact. For $1 \leq i \leq k$, choose $x_i \in A_i$. Then $\{ \mathcal{N}_\rho(x_i,n) : n \geq 1 \}$ is an open cover of $A_i$. Hence, there is an $n_i \geq 1$ such that $A_i \subset \mathcal{N}_\rho(x_i,n_i)$. Therefore, $\mathcal{N}_\rho(A_i,\epsilon) \subset \mathcal{N}_\rho(x_i,n_i + \epsilon)$. Consequently, $\mathrm{cl}(\mathcal{N}_\rho(A,\epsilon)) \subset \bigcup_{1 \leq i \leq k} \mathrm{cl}(\mathcal{N}_\rho(A_i,\epsilon)) \subset \bigcup_{1 \leq i \leq k} \mathrm{cl}(\mathcal{N}_\rho(x_i,n_i + \epsilon))$. It follows that $\mathrm{cl}(\mathcal{N}_\rho(A,\epsilon))$ is compact. ∎

**Definition.** Let $\mathcal{U}$ be an open cover of a topological space X and let $A \subset X$. Let $\overline{Star(A,\mathcal{U})} = \mathrm{cl}(\mathrm{Star}(A,\mathcal{U}))$, and for every $n \geq 1$, let $\overline{Star^n(A,\mathcal{U})} = \mathrm{cl}(\mathrm{Star}^n(A,\mathcal{U}))$.

**Definition.** An open cover $\mathcal{U}$ of a space X is *proper* if for every compact subset A of X, $\overline{\mathrm{Star}}(A,\mathcal{U})$ is compact.

**Lemma 3.2.** If $\mathcal{U}$ is a proper open cover of a topological space X, then for every compact subset A of X and for every $n \geq 1$ $\overline{\mathrm{Star}}^n(A,\mathcal{U})$ is compact.

**Proof.** We induct on n. Let $n \geq 1$ and assume $\overline{\mathrm{Star}}^n(A,\mathcal{U})$ is compact for every compact subset A of X. Let A be a compact subset of X. Since $\mathrm{Star}^{n+1}(A,\mathcal{U}) \subset \mathrm{Star}(\overline{\mathrm{Star}}^n(A,\mathcal{U}),\mathcal{U})$, then $\overline{\mathrm{Star}}^{n+1}(A,\mathcal{U}) \subset \overline{\mathrm{Star}}(\overline{\mathrm{Star}}^n(A,\mathcal{U}),\mathcal{U})$. Hence, $\overline{\mathrm{Star}}^{n+1}(A,\mathcal{U})$ is compact. ∎

**Definition.** A development $\{ \mathcal{U}_n : n \geq 1 \}$ in a topological space X is *proper* if $\mathcal{U}_n$ is proper for every $n \geq 1$.

Observe that if $\{ \mathcal{U}_n : n \geq 1 \}$ is a development in a topological space X in which $\mathcal{U}_{n+1}$ refines $\mathcal{U}_n$ for each $n \geq 1$, then $\{ \mathcal{U}_n : n \geq 1 \}$ is proper if and only if $\mathcal{U}_1$ is proper. In particular, if $\{ \mathcal{U}_n : n \geq 1 \}$ is a 2-development in a topological space X, then $\{ \mathcal{U}_n : n \geq 1 \}$ is proper if and only if $\mathcal{U}_1$ is proper.

**Example 3.1.** Let $\rho$ be a proper $\infty$-gauge on a topological space X. For every $n \geq 1$, let $\mathcal{U}_n = \{ \mathcal{N}_\rho(x,2^{-n}) : x \in X \}$. Then $\{ \mathcal{U}_n : n \geq 1 \}$ is a proper development in X. Indeed, $\mathcal{U}_{n+1}$ refines $\mathcal{U}_n$ for each $n \geq 1$. Also, for each $A \subset X$, $\mathrm{Star}(A,\mathcal{U}_1) \subset \mathcal{N}_\rho(A,1)$. Thus, if A is compact, then $\overline{\mathrm{Star}}(A,\mathcal{U}_1)$ is compact because $\overline{\mathrm{Star}}(A,\mathcal{U}_1) \subset \mathrm{cl}(\mathcal{N}_\rho(A,1))$ and $\mathrm{cl}(\mathcal{N}_\rho(A,1))$ is compact by Lemma 3.1.

**Theorem 3.3.** Let $\{ \mathcal{U}_n : n \geq 1 \}$ be a 3-development in a topological space X and let $\rho$ be the associated A-U distance function. Then $\{ \mathcal{U}_n : n \geq 1 \}$ is proper development if and only if $\rho$ is a proper $\infty$-gauge.

Our proof of Theorem 3.3 relies on:

**Lemma 3.4.** Suppose $\{ \mathcal{U}_n : n \geq 1 \}$ is a 3-development in a topological space X and $\rho$ is the associated A-U distance function. Then

*1)* for every subset A of X, $\mathrm{Star}(A,\mathcal{U}_1) \subset \mathcal{N}_\rho(A,1)$, and

*2)* for every $x \in X$ and every $n \geq 1$, $\mathcal{N}_\rho(x, {}^n/_2) \subset \text{Star}^{2n-1}(x, \mathcal{U}_1)$.

**Proof.** According to Theorem 2.1, $\rho$ has the property that $\text{Star}(x, \mathcal{U}_1) \subset \mathcal{N}_\rho(x, \frac{1}{2}) \subset \mathcal{N}_\rho(x, 1)$ for every $x \in X$. Hence, $\text{Star}(A, \mathcal{U}_1) \subset \mathcal{N}_\rho(A, 1)$ for every $A \subset X$. This proves assertion *1)*.

Assertion *2)* is proved by induction. Theorem 2.1 tells us that $\mathcal{N}_\rho(x, \frac{1}{2}) \subset \text{Star}^1(x, \mathcal{U}_1)$ for every $x \in X$. Let $n \geq 1$ and assume $\mathcal{N}_\rho(x, {}^n/_2) \subset \text{Star}^{2n-1}(x, \mathcal{U}_1)$. Let $y \in \mathcal{N}_\rho(x, {}^{(n+1)}/_2) - \mathcal{N}_\rho(x, {}^n/_2)$. Then, according to the recipe for $\rho$ formulated in Theorem 2.1, there is a $\mathcal{U}^*$-chain $C = (U_1, U_2, \cdots, U_k)$ joining $x$ to $y$ such that ${}^n/_2 \leq \lambda(C) < {}^{(n+1)}/_2$. It follows that there is a $j \in \{ 1, 2, \cdots, k \}$ such that $\Sigma_{1 \leq i \leq j-1} \delta(U_i) < {}^n/_2 \leq \Sigma_{1 \leq i \leq j} \delta(U_i)$. Hence, $\Sigma_{j+1 \leq i \leq k} \delta(U_i) < \frac{1}{2}$. Let $z \in U_{j-1} \cap U_j$ and $z' \in U_j \cap U_{j+1}$. Therefore, $\rho(x, z) < {}^n/_2$ and $\rho(z', y) < \frac{1}{2}$. Also, since $U_j$ is a subset of an element of $\mathcal{U}_1$, then $\{z, z'\} \subset V \in \mathcal{U}_1$. By the inductive hypothesis, $z \in \mathcal{N}_\rho(x, {}^n/_2) \subset \text{Star}^{2n-1}(x, \mathcal{U}_1)$. Therefore, $z' \in \text{Star}^{2n}(x, \mathcal{U}_1)$. Also, Theorem 2.1 tells us that $y \in \mathcal{N}_\rho(z', \frac{1}{2}) \subset \text{Star}^1(z', \mathcal{U}_1)$. Consequently, $y \in \text{Star}^{2n+1}(x, \mathcal{U}_1)$. This proves $\mathcal{N}_\rho(x, {}^{(n+1)}/_2) \subset \text{Star}^{2n+1}(x, \mathcal{U}_1)$. ∎

**Proof of Theorem 3.3.** First assume that $\{ \mathcal{U}_n : n \geq 1 \}$ is a proper development. Let $x \in X$ and $\epsilon > 0$. Let $n \geq 1$ such that $\epsilon \leq {}^n/_2$. Then Lemma 3.4 implies $\mathcal{N}_\rho(x, \epsilon) \subset \mathcal{N}_\rho(x, {}^n/_2) \subset \text{Star}^{2n-1}(x, \mathcal{U}_1)$. Lemma 3.2 implies that $\overline{\text{Star}^{2n-1}(x, \mathcal{U}_1)}$ is compact. Therefore, $\text{cl}(\mathcal{N}_\rho(x, \epsilon))$ is compact. This proves $\rho$ is a proper $\infty$-gauge.

Second assume $\rho$ is a proper $\infty$-gauge. Let $A$ be a compact subset of $X$. Lemma 3.4 implies that $\text{Star}(A, \mathcal{U}_1) \subset \mathcal{N}_\rho(A, 1)$. Lemma 3.1 tells us that $\text{cl}(\mathcal{N}_\rho(A, 1))$ is compact. Therefore, $\text{cl}(\text{Star}(A, \mathcal{U}_1))$ is compact. Hence, $\mathcal{U}_1$ is a proper open cover. It follows that $\{ \mathcal{U}_n : n \geq 1 \}$ is a proper development. ∎

**Corollary 3.5.** If $\mathbb{D}$ is a collection of proper 3-developments in a Hausdorff space $X$ that determines the given topology on $X$, then the collection of all the A-U distance functions associated with the elements of $\mathbb{D}$ is a proper $\infty$-gauge structure on $X$ that determines the given topology on $X$. ∎

To proceed beyond this point, we would like to transform proper $\infty$-gauges on the space $X$ into proper finite-valued gauges on $X$. We can't simply decapitate the proper $\infty$-gauges provided by Corollary 3.5, because if $X$ is non-compact, the decapitated gauges will not be proper. Instead, we will introduce another device called *tunnel systems* to obtain finite-valued gauges. Tunnel systems are a variant of the notion of *bridges* that is employed in [AMN].

**Definition.** Suppose $\rho$ is an $\infty$-gauge on a topological space $X$ and $\mathcal{Y}$ is the set of all $\rho$-crevasses. Let

$$\mathcal{Y}^{[2]} = \{ \{x, y\} : \text{there is a } Y \in \mathcal{Y} \text{ such that } x, y \in Y \}.$$

A *tunnel system for ρ* is an ordered pair $(\mathcal{T},\lambda)$ where $\mathcal{T} \subset \{ \{x,y\} : x,\ y \in X \}$ and $\lambda : \mathcal{T} \to (0,\infty)$ is a function such that the following three conditions are satisfied.

1) $\mathcal{T} \cap \mathcal{Y}^{[2]} = \emptyset$.

2) For all $x,\ y \in X$, there is a finite sequence $x = z_0,\ z_1,\ \cdots,\ z_k = y$ of points of X such that $\{z_{i-1},z_i\} \in \mathcal{Y}^{[2]} \cup \mathcal{T}$ for $1 \le i \le k$.

3) $\inf \lambda(\mathcal{T}) > 0$.

The elements of $\mathcal{T}$ are called *tunnels.*

Observe that if ρ is an ∞-gauge on a topological space X and $(\mathcal{T},\lambda)$ is a tunnel system for ρ, then for all $x,\ y \in X$, there is a finite sequence $x = z_0,\ z_1,\ \cdots,\ z_k = y$ of points of X such that for each $i \in \{ 1,\ 2,\ \cdots,\ k \}$, either $z_{i-1}$ and $z_i$ lie in the same ρ-crevasse or $\{z_{i-1},z_i\} \in \mathcal{T}$.

**Definition.** Suppose ρ is an ∞-gauge on a topological space X, $\mathcal{Y}$ is the set of all ρ-crevasses and $(\mathcal{T},\lambda)$ is a tunnel system for ρ. Define $\Lambda : \mathcal{Y}^{[2]} \cup \mathcal{T} \to [0,\infty)$ by

$$\Lambda(\{x,y\}) = \rho(x,y) \text{ if } \{x,y\} \in \mathcal{Y}^{[2]}, \text{ and } \Lambda(\{x,y\}) = \lambda(\{x,y\}) \text{ if } \{x,y\} \in \mathcal{T}.$$

Let $\mathcal{S} = \{ (z_0,\ z_1,\ \cdots,\ z_k) : k \ge 1 \text{ and } \{z_{i-1},z_i\} \in \mathcal{Y}^{[2]} \cup \mathcal{T} \text{ for } 1 \le i \le k \}$. Define $\mu : \mathcal{S} \to [0,\infty)$ by

$$\mu((z_0,\ z_1,\ \cdots,\ z_k)) = \Sigma_{1 \le i \le k} \Lambda(\{z_{i-1},z_i\}).$$

For $x,\ y \in X$, let $\mathcal{S}(x,y) = \{ (z_0,\ z_1,\ \cdots,\ z_k) \in \mathcal{S} : z_0 = x \text{ and } z_k = y \}$. Define *distance function* $\sigma : X \times X \to [0,\infty)$ *associated with* ρ and $(\mathcal{T},\lambda)$ by

$$\sigma(x,y) = \inf \{ \mu(S) : S \in \mathcal{S}(x,y) \}.$$

for $x,\ y \in X$.

Thus, $\sigma(x,y)$ is the infimum of all sums of lengths of sequences of "steps" joining x to y where each step either lies in a single ρ-crevasse or is a tunnel.

**Theorem 3.6.** If ρ is an ∞-gauge on a topological space X and $(\mathcal{T},\lambda)$ is a tunnel system for ρ, then the distance function σ associated with ρ and $(\mathcal{T},\lambda)$ is a (finite-valued) gauge on X that is equivalent to ρ such that $\sigma \le \rho$. Furthermore, if $\lambda_0 = \inf \lambda(\mathcal{T})$, then for all x, y $\in X$, $\sigma(x,y) = \rho(x,y)$ whenever either $\rho(x,y) < \lambda_0$ or $\sigma(x,y) < \lambda_0$.

**Proof.** It is obvious from the definition of σ that it satisfies conditions *1)*, *2)* and *3)* of the definition of a gauge.

To prove $\sigma \le \rho$, let $x,\ y \in X$. If $\{x,y\} \notin \mathcal{Y}^{[2]}$, then $\sigma(x,y) < \infty = \rho(x,y)$. If $\{x,y\} \in \mathcal{Y}^{[2]}$, then $(x,y) \in \mathcal{S}(x,y)$ and, hence, $\sigma(x,y) \le \mu((x,y)) = \Lambda(\{x,y\}) = \rho(x,y)$.

Next suppose $x \in X$ and $\epsilon \in (0,\infty)$. Let $y \in \mathcal{N}_\sigma(x,\epsilon)$. Since $\sigma \leq \rho$, then $\mathcal{N}_\rho(y,\epsilon - \sigma(x,y))$ $\subset \mathcal{N}_\sigma(y,\epsilon - \sigma(x,y)) \subset \mathcal{N}_\sigma(x,\epsilon)$. This proves $\mathcal{N}_\sigma(x,\epsilon)$ is an open subset of $X$. Therefore, $\sigma$ satisfies condition *4)* of the definition of a gauge.

Suppose $x, y \in X$ and either $\rho(x,y) < \lambda_0$ or $\sigma(x,y) < \lambda_0$. Since $\sigma \leq \rho$, then $\sigma(x,y) < \lambda_0$. Choose $\delta > 0$ so that $\sigma(x,y) + \delta < \lambda_0$. Then there is a finite sequence $(z_0, z_1, \cdots, z_k) \in \mathcal{S}(x,y)$ such that $\mu((z_0, z_1, \cdots, z_k)) < \sigma(x,y) + \delta$. Thus, $\Sigma_{1 \leq i \leq k} \Lambda(\{z_{i-1},z_i\}) < \sigma(x,y) + \delta$. Hence, $\Lambda(\{z_{i-1},z_i\}) < \lambda_0$ for $1 \leq i \leq k$. Consequently, $\Lambda(\{z_{i-1},z_i\}) = \rho(z_{i-1},z_i)$ for $1 \leq i \leq k$. Therefore, $\rho(x,y) = \rho(z_0,z_k) \leq \Sigma_{1 \leq i \leq k} \rho(z_{i-1},z_i) < \sigma(x,y) + \delta$. Since $\delta$ can be chosen to be arbitrarily small, then $\rho(x,y) \leq \sigma(x,y)$. This proves $\sigma(x,y) = \rho(x,y)$.

It follows that for $x \in X$ and $0 < \epsilon \leq \lambda_0$, $\mathcal{N}_\sigma(x,\epsilon) = \mathcal{N}_\rho(x,\epsilon)$. Therefore, $\rho$ and $\sigma$ are equivalent. ∎

**Example 3.2.** Let $\rho$ be an $\infty$-gauge on a topological space $X$ and let $\mathcal{Y}$ be the set of all $\rho$-crevasses. Choose $x_0 \in Y_0 \in \mathcal{Y}$, and for every $Y \in \mathcal{Y} - \{Y_0\}$, choose $x_Y \in Y$. Let $\mathcal{T} = \{\{x_0,x_Y\} : Y \in \mathcal{Y} - \{Y_0\}\}$, and define $\lambda : \mathcal{T} \to (0,\infty)$ by $\lambda(\{x_0,x_Y\}) = 1$ for each $Y \in \mathcal{Y} - \{Y_0\}$. Then $(\mathcal{T},\lambda)$ is a tunnel system for $\rho$. Hence, the distance function $\sigma$ associated with $\rho$ and $(\mathcal{T},\lambda)$ is a (finite-valued) gauge on $X$ that is equivalent to $\rho$. However, if $\mathcal{Y}$ is infinite, then $\sigma$ is not proper. Indeed, since $\mathrm{cl}(\mathcal{N}_\sigma(x_0,2)) \supset \{x_0\} \cup \{x_Y : Y \in \mathcal{Y} - \{Y_0\}\}$, then $\mathcal{Y}$ is an infinite open cover of $\mathrm{cl}(\mathcal{N}_\sigma(x_0,2))$ no proper subset of which covers $\mathrm{cl}(\mathcal{N}_\sigma(x_0,2))$.

**Definition.** Suppose $\rho$ is an $\infty$-gauge on a topological space $X$ and $(\mathcal{T},\lambda)$ is a tunnel system for $\rho$. For $A \subset X$ and $\epsilon \in (0,\infty)$, let

$T(A,\epsilon) = \{y \in X : \text{there is an } x \in A \text{ such that } \{x,y\} \in \mathcal{T} \text{ and } \lambda(\{x,y\}) < \epsilon\}$.

We say $(\mathcal{T},\lambda)$ is a *proper* tunnel system if $\mathrm{cl}(T(A,\epsilon))$ is compact for every compact subset $A$ of $X$ and every $\epsilon \in (0,\infty)$.

**Theorem 3.7.** Suppose $\rho$ is an $\infty$-gauge on a topological space $X$, $(\mathcal{T},\lambda)$ is a tunnel system for $\rho$, and $\sigma$ is the distance function associated with $\rho$ and $(\mathcal{T},\lambda)$. Then $\sigma$ is proper if and only if $\rho$ and $(\mathcal{T},\lambda)$ are proper.

**Proof.** First assume $\rho$ and $(\mathcal{T},\lambda)$ are proper. Let $\lambda_0 = \inf \lambda(\mathcal{T})$. Then $\lambda_0 > 0$. Let $x \in X$. We will prove by induction that for every $n \geq 1$, $\mathrm{cl}(\mathcal{N}_\sigma(x,n\lambda_0))$ is compact. To begin: since $\mathcal{N}_\sigma(x,\lambda_0) = \mathcal{N}_\rho(x,\lambda_0)$ and $\rho$ is proper, then $\mathrm{cl}(\mathcal{N}_\sigma(x,\lambda_0))$ is compact. Let $n \geq 1$ and assume that $\mathrm{cl}(\mathcal{N}_\sigma(x,n\lambda_0))$ is compact. Let $C = \mathrm{cl}(\mathcal{N}_\sigma(x,n\lambda_0))$. We assert that

$$\mathcal{N}_\sigma(x,(n+1)\lambda_0) \subset \mathcal{N}_\rho(C,(n+2)\lambda_0) \cup \mathcal{N}_\rho(T(C,(n+1)\lambda_0),\lambda_0).$$

To prove this assertion, let $y \in \mathcal{N}_\sigma(x,(n+1)\lambda_0)$. If $\sigma(x,y) < n\lambda_0$, then $y \in C \subset \mathcal{N}_\rho(C,(n+2)\lambda_0)$. So we can assume $\sigma(x,y) \geq n\lambda_0$. Since $\sigma(x,y) < (n+1)\lambda_0$, then (using the terminology introduced in the definition of $\sigma$) there is a $(z_0, z_1, \cdots, z_k) \in \mathcal{S}(x,y)$ such that $n\lambda_0 \leq \mu((z_0, z_1, \cdots, z_k)) < (n+1)\lambda_0$. Since $\mu((z_0, z_1, \cdots, z_k)) = \Sigma_{1 \leq i \leq k} \Lambda(\{z_{i-1},z_i\})$, then there is a $j \in \{1, 2, \cdots, k\}$ such that $\Sigma_{1 \le i \le j-1} \Lambda(\{z_{i-1}, z_i\}) < n\lambda_0 \le \Sigma_{1 \le i \le j} \Lambda(\{z_{i-1}, z_i\}) \le \Sigma_{1 \le i \le k} \Lambda(\{z_{i-1}, z_i\}) < (n+1)\lambda_0.$

Hence,

$$\Sigma_{j+1 \le i \le k} \Lambda(\{z_{i-1}, z_i\}) < \lambda_0.$$

Since $\sigma(x, z_{j-1}) \le \Sigma_{1 \le i \le j-1} \Lambda(\{z_{i-1}, z_i\}) < n\lambda_0$, then $z_{j-1} \in C$. Since $\sigma(z_j, y) \le \Sigma_{j+1 \le i \le k} \Lambda(\{z_{i-1}, z_i\}) < \lambda_0$, then $\rho(z_j, y) = \sigma(z_j, y) < \lambda_0$. We now consider two cases.

**Case 1:** $\{z_{j-1}, z_j\} \in \mathcal{Y}^{[2]}$. In this case, $\rho(z_{j-1}, z_j) = \Lambda(\{z_{j-1}, z_j\}) \le \mu((z_0, z_1, \cdots, z_k)) < (n+1)\lambda_0$. Then $\rho(z_{j-1}, y) \le \rho(z_{j-1}, z_j) + \rho(z_j, y) < (n+2)\lambda_0$. Therefore, $y \in \mathcal{N}_\rho(C, (n+2)\lambda_0)$.

**Case 2:** $\{z_{j-1}, z_j\} \in \mathcal{T}$. In this case, $\lambda(z_{j-1}, z_j) = \Lambda(\{z_{j-1}, z_j\}) \le \mu((z_0, z_1, \cdots, z_k)) < (n+1)\lambda_0$. Hence, $z_j \in T(C, (n+1)\lambda_0)$. Therefore, $y \in \mathcal{N}_\rho(T(C, (n+1)\lambda_0), \lambda_0)$.

Thus, in all possible situations, $y \in \mathcal{N}_\rho(C, (n+2)\lambda_0) \cup \mathcal{N}_\rho(T(C, (n+1)\lambda_0), \lambda_0)$, proving our assertion.

Since C is compact and $\rho$ is proper, then Lemma 3.1 implies that $\mathcal{N}_\rho(C, (n+2)\lambda_0)$ has compact closure. Since $(\mathcal{T}, \lambda)$ is a proper tunnel system, then cl$(T(C, (n+1)\lambda_0))$ is compact. Since $\rho$ is proper, then Lemma 3.1 implies $\mathcal{N}_\rho(\text{cl}(T(C, (n+1)\lambda_0)), \lambda_0)$ has compact closure. Consequently, $\mathcal{N}_\rho(T(C, (n+1)\lambda_0), \lambda_0)$ has compact closure. Our assertion now implies that $\mathcal{N}_\sigma(x, (n+1)\lambda_0)$ has compact closure. It follows by induction that $\sigma$ is proper.

Now assume that $\sigma$ is proper. Since $\sigma \le \rho$, then $\mathcal{N}_\rho(x, \epsilon) \subset \mathcal{N}_\sigma(x, \epsilon)$ for each $x \in X$ and each $\epsilon > 0$. Therefore, the properness of $\sigma$ implies the properness of $\rho$. Let A be a compact subset of X and let $\epsilon > 0$. If $y \in T(A, \epsilon)$, then there is an $x \in A$ such that $\{x, y\} \in \mathcal{T}$ and $\lambda(\{x, y\}) < \epsilon$. Since $\sigma(x, y) \le \mu((x, y)) = \Lambda(\{x, y\}) = \lambda(\{x, y\})$, then $\sigma(x, y) < \epsilon$. This proves $T(A, \epsilon) \subset \mathcal{N}_\sigma(A, \epsilon)$. Since $\sigma$ is proper, then Lemma 3.1 implies $\mathcal{N}_\sigma(A, \epsilon)$ has compact closure. Consequently, $T(A, \epsilon)$ has compact closure. Therefore, $(\mathcal{T}, \lambda)$ is proper. ∎

We close this section by giving two examples of proper tunnel systems for an ∞-gauge $\rho$. Both these examples require the set of all $\rho$-crevasses to be countable. This is a necessary condition if we intend to find a proper (finite-valued) gauge that is equivalent to $\rho$, as the following proposition reveals.

**Proposition 3.8.** Let $\rho$ be an ∞-gauge on a topological space X. If X is $\sigma$-compact, then the set of all $\rho$-crevasses is countable. Therefore, if X has a proper (finite-valued) gauge, then the set of all $\rho$-crevasses is countable.

**Proof.** Suppose X is a $\sigma$-compact space. Then $X = \cup_{n \ge 1} C_n$ where each $C_n$ is compact. For each $n \ge 1$, let $\mathcal{Y}_n = \{ Y \in \mathcal{Y} : Y \cap C_n \ne \emptyset \}$. Then each $\mathcal{Y}_n$ is finite and $\mathcal{Y} = \cup_{n \ge 1} \mathcal{Y}_n$. Hence, $\mathcal{Y}$ is countable.

If X has a proper (finite-valued) gauge, then, as we observe previously, X is σ-compact. Therefore, we conclude that the set of all ρ-crevasses is countable. ∎

**Example 3.3.** Let ρ be an ∞-gauge on a Hausdorff space X, and assume that the set $\mathcal{Y}$ of all ρ-crevasses is countable. Say $\mathcal{Y} = \{ Y_1, Y_2, Y_3, \cdots \}$. For each $i \geq 1$, choose $x_i \in Y_i$. Let $\mathcal{T} = \{ \{x_i, x_{i+1}\} : i \geq 1 \}$ and define $\lambda : \mathcal{T} \rightarrow (0, \infty)$ by $\lambda(\{x_i, x_{i+1}\}) = 1$ for $i \geq 1$. Then $(\mathcal{T}, \lambda)$ is a proper tunnel system for ρ. Indeed, if A is a compact subset of X and $\epsilon > 0$, then there is an $n \geq 1$ such that $A \subset \cup_{1 \leq i \leq n} Y_i$. Then $T(A, \epsilon) \subset \{ x_1, x_2, \cdots, x_{n+1} \}$. Thus, $T(A, \epsilon)$ is a finite and, hence, a closed and compact subset of X.

**Example 3.4.** Let ρ be an ∞-gauge on a Hausdorff space X, and assume that the set $\mathcal{Y}$ of all ρ-crevasses is countable. Say $\mathcal{Y} = \{ Y_0, Y_1, Y_2, \cdots \}$. For each $i \geq 0$, choose $x_i \in Y_i$. Let $\mathcal{T} = \{ \{x_0, x_i\} : i \geq 1 \}$ and define $\lambda : \mathcal{T} \rightarrow (0, \infty)$ by $\lambda(\{x_0, x_i\}) = i$ for $i \geq 1$. Then $(\mathcal{T}, \lambda)$ is a proper tunnel system for ρ. Indeed, if $A \subset X$ and $\epsilon > 0$, then there is an $n \geq 0$ such that $n < \epsilon \leq n + 1$. Then $T(A, \epsilon) \subset \{ x_0, x_1, \cdots, x_n \}$. Thus, $T(A, \epsilon)$ is a finite and, hence, a closed and compact subset of X.

## 4. Invariant gauge and proper gauge structures

**Definition.** Let G be a group of homeomorphisms of a topological space X. An ∞-gauge ρ on X is *G-invariant* if $\rho(g(x), g(y)) = \rho(x, y)$ for all $g \in G$ and all x, y ∈ X. An open cover $\mathcal{U}$ of X is *G-invariant* if for every U ∈ $\mathcal{U}$ and every g ∈ G, g(U) ∈ $\mathcal{U}$. A development $\{ \mathcal{U}_n : n \geq 1 \}$ in X is *G-invariant* if $\mathcal{U}_n$ is G-invariant for every $n \geq 1$.

**Lemma 4.1.** If G is a group of homeomorphisms of a topological space X and $\{ \mathcal{U}_n : n \geq 1 \}$ is a G-invariant 3-development in X, then the associated A-U distance function and its decapitation are G-invariant.

**Proof.** This lemma follows immediately from the observation that if $\{ \mathcal{U}_n : n \geq 1 \}$ is G-invariant, then the recipe for the A-U distance function associated with $\{ \mathcal{U}_n : n \geq 1 \}$ yields an ∞-gauge on X which is clearly G-invariant, and the formula for its decapitation then yields a gauge on X that is clearly G-invariant gauge. ∎

**Definition.** Let G be a group of homeomorphisms of a topological space X, let ρ be an ∞-metric on X, and let $(\mathcal{T}, \lambda)$ be a tunnel system for ρ. $(\mathcal{T}, \lambda)$ is *G-invariant* if for all $\{x, y\} \in \mathcal{T}$ and all g ∈ G, $\{g(x), g(y)\} \in \mathcal{T}$ and $\lambda(\{g(x), g(y)\}) = \lambda(\{x, y\})$.

**Lemma 4.2.** If G is a group of homeomorphisms of a topological space X, ρ is a G-invariant ∞-metric on X, and $(\mathcal{T}, \lambda)$ is a G-invariant tunnel system for ρ, then the associated distance function σ is G-invariant.

**Proof.** It is clear that if ρ and $(\mathcal{T},\lambda)$ are G-invariant, then the recipe for creating σ from ρ and $(\mathcal{T},\lambda)$ yields a G-invariant gauge on X. ∎

**Notation.** If G is a group of homeomorphisms of a topological space X and A ⊂ X, let

$$GA = \textstyle\bigcup_{g \in G} g(A).$$

Recall that a group G of homeomorphisms of a space X is *nearly proper* if for all compact subsets A and B of X, cl($\bigcup$ { g(A) : g ∈ G and g(A) ∩ B ≠ ∅ }) is compact.

**Lemma 4.3.** Suppose G is a nearly proper group of homeomorphisms of a topological space X, ρ is a G-invariant ∞-metric on X, $\mathcal{Y}$ is the set of all ρ-crevasses, and $(\mathcal{T},\lambda)$ is a tunnel system for ρ with following property.

   For every compact subset A of X and for every ε > 0, cl(T(GA,ε)) is compact.

Let G$\mathcal{T}$ = { {g(x),g(y)} : g ∈ G and {x,y} ∈ $\mathcal{T}$ } and define $\lambda_G$ : G$\mathcal{T}$ → (0,∞) by

   $\lambda_G$({x,y})  =  inf { λ({x′,y′}) : {x′,y′} ∈ $\mathcal{T}$ and {x,y} = {g(x′),g(y′)} for some g ∈ G }.

Then (G$\mathcal{T}$,$\lambda_G$) is a G-invariant proper tunnel system for ρ.

**Proof.** Recall that $\mathcal{Y}^{[2]}$ = { {x,y} : there is a Y ∈ $\mathcal{Y}$ such that x, y ∈ Y }. Since ρ is G-invariant, then the elements of G permute the elements of $\mathcal{Y}$. Hence, for g ∈ G and x, y ∈ X, {x,y} ∈ $\mathcal{Y}^{[2]}$ if and only if {g(x),g(y)} ∈ $\mathcal{Y}^{[2]}$. Therefore, since $\mathcal{Y}^{[2]}$ ∩ $\mathcal{T}$ = ∅, then $\mathcal{Y}^{[2]}$ ∩ G$\mathcal{T}$ = ∅. Hence, (G$\mathcal{T}$,$\lambda_G$) satisfies condition *1*) of the definition of "tunnel system".

(G$\mathcal{T}$,$\lambda_G$) satisfies conditions *2*) and *3*) of the definition of "tunnel system" because $(\mathcal{T},\lambda)$ satisfies these conditions. Also, it is clear that (G$\mathcal{T}$,$\lambda_G$) is G-invariant.

It remains to verify that (G$\mathcal{T}$,$\lambda_G$) is proper. To this end, for A ⊂ X and ε > 0, let

   $T_G$(A,ε)  = { y ∈ X : there is an x ∈ A such that {x,y} ∈ G$\mathcal{T}$ and $\lambda_G$({x,y}) < ε }.

We must verify that if A is a compact subset of X and ε > 0, then cl($T_G$(A,ε)) is compact. Observe that for A ⊂ X and ε > 0,

$$T_G(A,\varepsilon) \;=\; \textstyle\bigcup_{g \in G} g(T(g^{-1}(A),\varepsilon)).$$

Also observe that for S ⊂ X and ε > 0, if T(S,ε) ≠ ∅, then S ∩ T(T(S,ε),ε) ≠ ∅. Suppose A is a compact subset of X and ε > 0. Then, by hypothesis, the sets B = cl(T(GA,ε)) and C = B ∪ cl(T(GB,ε)) are compact. Let g ∈ G and assume that T(g⁻¹(A),ε) ≠ ∅. Then g⁻¹(A) ∩ T(T(g⁻¹(A),ε),ε) ≠ ∅. Also T(g⁻¹(A),ε) ⊂ B. Hence, T(T(g⁻¹(A),ε),ε) ⊂ C. Thus, g⁻¹(A) ∩ C ≠ ∅. Therefore, g(T(g⁻¹(A),ε)) ⊂ g(B) ⊂ g(C) and A ∩ g(C) ≠ ∅. It follows that $T_G$(A,ε) ⊂ $\bigcup$ { g(C) : g ∈ G and g(C) ∩ A ≠ ∅ }. Since G is nearly proper, we conclude that cl($T_G$(A,ε)) is compact. ∎

**Corollary 4.4.** If G is a nearly proper group of homeomorphisms of a Hausdorff space X, and ρ is a G-invariant ∞-metric on X such that the set of all ρ-crevasses is countable, then there is a tunnel system $(\mathcal{T},\lambda)$ for ρ such that $(G\mathcal{T},\lambda_G)$ is a G-invariant proper tunnel system for ρ.

**Proof.** Let $(\mathcal{T},\lambda)$ be the tunnel system for ρ defined in Example 3.4. Then for every subset A of X and every ε > 0, T(GA,ε) is a finite set. Hence, cl(T(GA,ε)) = T(GA,ε) is compact. Therefore, Lemma 4.3 implies $(G\mathcal{T},\lambda_G)$ is a G-invariant and proper. ∎

## 5. Proofs of the Isometrization Theorems

Our proofs of the isometrization theorems depend on two key lemmas which allow us to apply the methods developed in sections 2, 3 and 4.

Recall that a group of homeomorphisms G of a topological space X is *equiregular* if for every x ∈ X and every open neighborhood U of x in X there is an open neighborhood V of x in X such that cl(V) ⊂ U and every y ∈ X has an open neighborhood $N_y$ with the property that for every g ∈ G, if $g(N_y) \cap$ cl(V) ≠ ∅, then $g(N_y) \subset$ U.

Here is the first of the two key lemmas.

**Lemma 5.1.** If G is an equiregular group of homeomorphisms of a topological space X and G\X is a paracompact regular space, then every G-invariant open cover of X is star-refined by a G-invariant open cover of X.

**Proof.** Let $\mathcal{U}$ be a G-invariant open cover of X. Since G is equiregular, there is an open cover $\mathcal{V}$ of X with the property that for every V ∈ $\mathcal{V}$, there is a $U_V \in \mathcal{U}$ such that cl(V) ⊂ $U_V$ and for every x ∈ X, there is an open neighborhood N of x in X such that for every g ∈ G, if g(N) ∩ cl(V) ≠ ∅, then g(N) ⊂ $U_V$. (The choice of N depends on V, $U_V$ and x.) Since π : X → G\X is an open map and G\X is paracompact and regular, then there is a locally finite open cover $\mathcal{L}$ of G\X such that { cl(L) : L ∈ $\mathcal{L}$ } refines { π(V) : V ∈ $\mathcal{V}$ }. For each L ∈ $\mathcal{L}$, choose V(L) ∈ $\mathcal{V}$ such that cl(L) ⊂ π(V(L)). For each x ∈ X, let $\mathcal{L}$(x) = { L ∈ $\mathcal{L}$ : π(x) ∈ cl(L) } and let C(x) = $\cup_{L \in \mathcal{L} - \mathcal{L}(x)}$cl(L). Then for each x ∈ X, $\mathcal{L}$(x) is a finite subset of $\mathcal{L}$ and C(x) is a closed subset of G\X such that π(x) ∉ C(x).

For every x ∈ X and every L ∈ $\mathcal{L}$(x), since π(x) ∈ cl(L) ⊂ π(V(L)), then there is a $g_{x,L} \in$ G such that x ∈ $g_{x,L}$(V(L)). Since G is equiregular, then for each x ∈ X and each L ∈ $\mathcal{L}$(x), there is an open neighborhood $N_{x,L}$ of x in X such that $N_{x,L} \subset g_{x,L}$(V(L)) and for every g ∈ G, if $g(N_{x,L}) \cap$ cl(V(L)) ≠ ∅, then $g(N_{x,L}) \subset U_{V(L)}$. Let $N_x = (\cap_{L \in \mathcal{L}(x)}N_{x,L}) - \pi^{-1}$(C(x)). Then $N_x$ is an open neighborhood of x in X. Furthermore, for x ∈ X and g ∈ G: observe that if L ∈ $\mathcal{L}$(x) and $g(N_x) \cap$ cl(V(L)) ≠ ∅, then $g(N_x) \subset U_{V(L)}$; and observe that if L ∈ $\mathcal{L}$ and $g(N_x) \cap \pi^{-1}$(cl(L)) ≠ ∅, then L ∉ C(x) and, hence, L ∈ $\mathcal{L}$(x).

Let $\mathcal{W} = \{ g(N_x) : x \in X$ and $g \in G \}$. The $\mathcal{W}$ is clearly a G-invariant open cover of X. It remains to prove that $\mathcal{W}$ star-refines $\mathcal{U}$. Let $x \in X$. Choose $L \in \mathcal{L}(x)$. We will prove that there is an $h \in G$ such that $\text{Star}(x,\mathcal{W}) \subset h(U_{V(L)})$. Since $\mathcal{U}$ is G-invariant, it will then follow that $\mathcal{W}$ star-refines $\mathcal{U}$. Choose a specific $y \in X$ and a specific $g \in G$ such that $x \in g(N_y)$. Then $x \in g(N_y) \cap \pi^{-1}(\text{cl}(L))$. Therefore, $\pi(N_y) \cap \text{cl}(L) \neq \emptyset$. Consequently, $L \in \mathcal{L}(y)$. Thus, $N_y \subset N_{y,L} \subset g_{y,L}(V(L))$. Therefore, $x \in g(N_y) \subset gg_{y,L}(V(L))$. Let $h = gg_{y,L}$. Then $h \in G$ and $x \in h(V(L))$. Now suppose $y' \in X$ and $g' \in G$ are arbitrary elements of X and G, respectively, such that $x \in g'(N_{y'})$. Then $x \in g'(N_{y'}) \cap h(V(L))$. Therefore, $h^{-1}g'(N_{y'}) \cap V(L) \neq \emptyset$. It follows that $h^{-1}g'(N_{y'}) \subset U_{V(L)}$. Hence, $g'(N_{y'}) \subset h(U_{V(L)})$. This proves $\text{Star}(x,\mathcal{W}) \subset h(U_{V(L)})$. ∎

**Remark.** This proof is takes inspiration from A. H. Stone's proof [S] that in a paracompact regular space, every open cover is star-refined by an open cover.

We restate:

**The Isometrization Theorem 1.1.** If X is a Hausdorff space, G is a group of homeomorphisms of X, and G\X is a paracompact regular space, then: G is isometrizable if and only if G is equiregular.

**Proof.** First assume G is isometrizable. Let $\mathcal{P}$ be a G-invariant gauge structure on X. As we observed in section 1, we can assume that $\mathcal{P}$ is closed under maximization. Let $x \in X$ and let U be an open neighborhood of x in X. Then there is a $\rho \in \mathcal{P}$ and an $\epsilon \in (0,\infty)$ such that $\mathcal{N}_\rho(x,\epsilon) \subset U$. Let $V = \mathcal{N}_\rho(x,\epsilon/3)$. Then $\text{cl}(V) \subset \mathcal{N}_\rho(x,\epsilon/3] \subset U$. Also for each $y \in Y$, let $N_y = \mathcal{N}_\rho(y,\epsilon/3)$. Let $g \in G$ and assume $g(N_y) \cap \text{cl}(V) \neq \emptyset$. Since $\rho$ is G-invariant, then $g(N_y) = \mathcal{N}_\rho(g(y),\epsilon/3)$. Thus, $\emptyset \neq g(N_y) \cap \text{cl}(V) \subset \mathcal{N}_\rho(g(y),\epsilon/3) \cap \mathcal{N}_\rho(x,\epsilon/3]$. Consequently, $g(N_y) \subset \mathcal{N}_\rho(x,\epsilon) \subset U$. This proves G is equiregular.

Second assume G is equiregular. Let $x \in X$ and let U be an open neighborhood of x in X. Then there is an open neighborhood V of x in X such that $\text{cl}(V) \subset U$ and each point $y \in X$ has an open neighborhood $N_y$ such that for every $g \in G$, if $g(N_y) \cap \text{cl}(V) \neq \emptyset$, then $g(N_y) \subset U$. Let $\mathcal{V}_{(x,U,1)} = \{ g(N_y) : y \in X$ and $g \in G \}$. Then Let $\mathcal{V}_{(x,U,1)}$ is a G-invariant open cover of X such that $\text{Star}(x,\mathcal{V}_{(x,U,1)}) \subset U$. Repeated application of Lemma 5.1 provides a sequence $\mathcal{V}_{(x,U,2)}, \mathcal{V}_{(x,U,3)}, \mathcal{V}_{(x,U,4)}, \cdots$ of G-invariant open covers of X such that for each $i \geq 1$, $\mathcal{V}_{(x,U,i+1)}$ star-refines $\mathcal{V}_{(x,U,i)}$. Hence, $\mathcal{D}_{(x,U)} = \{ \mathcal{V}_{(x,U,1)}, \mathcal{V}_{(x,U,3)}, \mathcal{V}_{(x,U,5)}, \cdots \}$ is a G-invariant 3-development in X such that $\text{Star}(x,\mathcal{V}_{(x,U,1)}) \subset U$. Let $\sigma_{(x,U)}$ be the decapitation of the A-U distance function associated with $\mathcal{D}_{(x,U)}$. Then Corollary 2.2 and Lemma 4.1 imply that $\sigma_{(x,U)}$ is a G-invariant (finite-valued) gauge on X; and Theorem 2.1 implies $\mathcal{N}_{\sigma_{(x,U)}}(x,\frac{1}{2}) \subset \text{Star}(x,\mathcal{V}_{(x,U,1)}) \subset U$. Thus, $\{ \sigma_{(x,U)} : x \in X$ and U is an open neighborhood of x in X $\}$ is a G-invariant gauge structure on X. This proves G is isometrizable. ∎

Before proving the second of the two key lemmas, we need the following preliminary result. We alert the reader to the fact that in a non-Hausdorff space, a compact subset may not be closed.  We remind the reader that a space is locally compact if every point has an open neighborhood whose closure is compact.  A space is σ-compact if it is the union of a countable collection of compact subsets.

**Proposition 5.2.**  Suppose Y is a σ-compact regular space with the property that every point lies in the interior of a compact subset.  Then Y is locally compact and paracompact.  Moreover, there is a sequence { $K_i$ } of closed compact subsets of Y such that Y = $\cup_{i \geq 1} K_i$ and $K_i \cap K_j = \emptyset$ whenever | i − j | > 1.  Also, there is a sequence { $L_i$ } of open subsets of Y such that $K_i \subset L_i$, $L_i \cap K_j = \emptyset$ whenever | i − j | > 1, and $L_i \cap L_j = \emptyset$ whenever | i − j | > 2.

**Proof.**  Since Y is regular, it follows that every point of Y has an open neighborhood with compact closure.  In other words, Y is locally compact.  Hence, every compact subset of Y lies in the interior of a closed compact set.  Since Y is σ-compact, there is a sequence { $C_i$ } of (possibly non-closed) compact subsets of Y such that Y = $\cup_{i \geq 1} C_i$. Inductively choose a sequence { $D_i$ } of closed compact subsets of Y such that $C_1 \subset$ int($D_1$) and $C_i \cup D_{i-1} \subset$ int($D_i$) for each i ≥ 2.  Then Y = $\cup_{i \geq 1} D_i$ and $D_i \subset$ int($D_{i+1}$) for each i ≥ 1.

Let $K_1 = D_1$ and for each i ≥ 2, let $K_i = D_i -$ int($D_{i-1}$).  Also let $L_1 =$ int($D_2$), $L_2 =$ int($D_3$) and for each i ≥ 3, let $L_i =$ int($D_{i+1}$) − $D_{i-2}$.  Then { $K_i$ } is a sequence of closed compact subsets of Y such that Y = $\cup_{i \geq 1} K_i$ and $K_i \cap K_j = \emptyset$ whenever | i − j | > 1.  Also { $L_i$ } is a sequence of open subsets of Y such that $K_i \subset L_i$, $L_i \cap K_j = \emptyset$ whenever | i − j | > 1, and $L_i \cap L_j = \emptyset$ whenever | i − j | > 2.

To prove that Y is paracompact, let $\mathcal{U}$ be an open cover of Y.  For each i ≥ 1, let $\mathcal{U}_i$ be a finite subset of $\mathcal{U}$ that covers $K_i$ and let $\mathcal{V}_i =$ { U ∩ $L_i$ : U ∈ $\mathcal{U}_i$ }.  Let $\mathcal{V} = \mathcal{V}_1 \cup \mathcal{V}_2 \cup \mathcal{V}_3 \cup$ ⋯.  Then $\mathcal{V}$ is an open cover of Y that refines $\mathcal{U}$.  $\mathcal{V}$ is locally finite because { $L_i$ } is an open cover of Y and { V ∈ $\mathcal{V}$ : $L_i \cap$ V ≠ $\emptyset$ } is a subset of the finite set $\mathcal{V}_{i-2} \cup \mathcal{V}_{i-1} \cup \mathcal{V}_i \cup \mathcal{V}_{i+1} \cup \mathcal{V}_{i+2}$. ∎

We now state the second key lemma.

**Lemma 5.3.**  If X is a locally compact σ-compact Hausdorff space, G is a nearly proper group of homeomorphisms of X, and G\X is a regular space, then X has a G-invariant proper open cover.

**Proof.**  Since the quotient map π : X → G\X is a continuous open map, then G\X is σ-compact and every point of G\X lies in the interior of a compact set.  Thus, Proposition 5.2 implies that G\X is locally compact and paracompact, and there is a sequence { $K_i$ } of closed compact subsets of G\X such that G\X = $\cup_{i \geq 1} K_i$ and $K_i \cap K_j = \emptyset$ whenever | i − j | > 1, and there is a sequence { $L_i$ } of open subsets of G\X such that $K_i \subset L_i$, $L_i \cap K_j = \emptyset$ whenever | i − j | > 1, and  $L_i \cap L_j = \emptyset$ whenever | i − j | > 2.

For every i ≥ 1, there is a collection $\mathcal{U}_i$ of open subsets of X with compact closure such that $K_i \subset \cup \{ \pi(U) : U \in \mathcal{U}_i \} \subset L_i$. Since $K_i$ is compact, then there is a finite subset $\mathcal{F}_i$ of $\mathcal{U}_i$ such that $\cup \{ \pi(U) : U \in \mathcal{F}_i \}$ covers $K_i$. Let $\mathcal{F} = \mathcal{F}_1 \cup \mathcal{F}_2 \cup \mathcal{F}_3 \cup \cdots$. Then $\{ \pi(U) : U \in \mathcal{F} \}$ covers G\X. Let $\mathcal{V} = \{ g(U) : U \in \mathcal{F}$ and $g \in G \}$. Then $\mathcal{V}$ is a G-invariant open cover of X.

To prove that $\mathcal{V}$ is proper, let B be a compact subset of X. Then there is a k ≥ 1 such that $\pi(B) \cap L_i = \emptyset$ for all i > k. Hence, if i > k and $U \in \mathcal{F}_i$, then $\pi(B) \cap \pi(U) = \emptyset$ because $\pi(U) \subset L_i$. Therefore, if i ≥ 1, $U \in \mathcal{F}_i$, $g \in G$ and $B \cap g(U) \neq \emptyset$, then i ≤ k. Thus,

$$\{ U \in \mathcal{F} : B \cap g(U) \neq \emptyset \text{ for some } g \in G \} \subset \mathcal{F}_1 \cup \mathcal{F}_2 \cup \cdots \cup \mathcal{F}_k.$$

Let $A = \cup \{ cl(U) : U \in \mathcal{F}_1 \cup \mathcal{F}_2 \cup \cdots \cup \mathcal{F}_k \}$. Then A is compact because each $\mathcal{F}_i$ is finite and its elements have compact closure. Also if $U \in \mathcal{F}$, $g \in G$ and $B \cap g(U) \neq \emptyset$, then $U \subset A$. Thus, if $U \in \mathcal{F}$, $g \in G$ and $B \cap g(U) \neq \emptyset$, then $g(U) \subset g(A)$. Therefore, $Star(B,\mathcal{V}) \subset \cup \{ g(A) : g \in G$ and $g(A) \cap B \neq \emptyset \}$. Since G is nearly proper, then $cl(\cup \{ g(A) : g \in G$ and $g(A) \cap B \neq \emptyset \})$ is compact. Hence, $cl(Star(B,\mathcal{V}))$. This proves $\mathcal{V}$ is proper. ∎

We restate:

**The Proper Isometrization Theorem 1.3.** If X is a locally compact σ-compact Hausdorff space and G\X is a regular space, then: G is properly isometrizable if and only if G is equiregular and nearly proper.

**Proof.** First assume G is properly isometrizable. Let $\mathcal{P}$ be a G-invariant proper gauge structure on X. The Isometrization Theorem implies that G is equiregular. To prove that G is nearly proper, let A and B be compact subsets of X. Let $\rho \in \mathcal{P}$ and let $x_0 \in X$. Since $\{ \mathcal{N}_\rho(x_0,n) : n \geq 1 \}$ is an increasing open cover of X, then there is an n ≥ 1 such that $A \cup B \subset \mathcal{N}_\rho(x_0,n)$. Since $\rho$ is G-invariant, then for every $g \in G$, $g(A) \subset \mathcal{N}_\rho(g(x_0),n)$. Hence,

$\cup \{ g(A) : g \in G$ and $g(A) \cap B \neq \emptyset \}$

$\subset \cup \{ \mathcal{N}_\rho(g(x_0),n) : g \in G$ and $\mathcal{N}_\rho(g(x_0),n) \cap \mathcal{N}_\rho(x_0,n) \neq \emptyset \}$

$\subset \mathcal{N}_\rho(x_0,3n).$

Since $\rho$ is proper, then $cl(\mathcal{N}_\rho(x_0,3n))$ is compact. Consequently, $cl(\cup \{ g(A) : g \in G$ and $g(A) \cap B \neq \emptyset \})$ is compact. This proves G is nearly proper.

Second assume G is equiregular and nearly proper. Since the quotient map $\pi : X \to$ G\X is continuous and open, then G\X is σ-compact and every point of G\X lies in the interior of a compact set. Thus, Proposition 5.2 implies that G\X is paracompact. Hence the Isometrization Theorem provides a G-invariant gauge structure $\mathcal{P}$.

Lemma 5.3 provides a G-invariant proper open cover $\mathcal{U}_1$ of X. Repeated application of Lemma 5.1 provides a sequence $\mathcal{U}_2, \mathcal{U}_3, \mathcal{U}_4, \cdots$ of G-invariant open covers of X such that for each i ≥ 1, $\mathcal{U}_{i+1}$ star-refines $\mathcal{U}_i$. Thus, { $\mathcal{U}_{2n-1}$ : n ≥ 1 } is a G-invariant proper 3-development in X. Let σ be the A-U distance function associated with { $\mathcal{U}_{2n-1}$ : n ≥ 1 }. Then Theorem 3.3 and Lemma 4.1 imply that σ is a G-invariant proper ∞-gauge on X. Since X is σ-compact, then Proposition 3.8 tells us that the set of all σ-crevasses is countable. Since G is nearly proper, then Corollary 4.4 provides a G-invariant proper tunnel system (G$\mathcal{T}$,λ$_G$) for σ. Let τ be the distance function associated with σ and (G$\mathcal{T}$,λ$_G$). Theorem 3.6 implies that τ is a (finite-valued) gauge on X that is equivalent to σ. Since σ and (G$\mathcal{T}$,λ$_G$) are G-invariant and proper, then Theorem 3.7 and Lemma 4.2 imply that τ is a G-invariant and proper. Observations from section 1 tell us that since τ is proper, then for every ρ ∈ $\mathcal{P}$, max{ρ,τ} is a G-invariant proper gauge on X that is equivalent to { ρ, τ }. Hence, { max{ρ,τ} : ρ ∈ $\mathcal{P}$ } is a G-invariant proper gauge structure on X. This proves G is properly isometrizable. ∎

## 6. Groups of homeomorphisms that act properly

Our goal in this section is to prove that Theorem 1.5 [AMN] and Corollary 1.6 [AMN] can be deduced from The Proper Isometrization Theorem 1.3 and Corollary 1.4. We restate the theorems from [AMN].

**Theorem 1.5. [AMN].** If G is a group of homeomorphisms of a locally compact σ-compact Hausdorff space X and G acts properly on X, then G is properly isometrizable.

**Corollary 1.6. [AMN].** If G is a group of homeomorphisms of a locally compact σ-compact metrizable space X and G acts properly on X, then G is singly properly isometrizable.

Clearly, to achieve this goal, it suffices to show that if G is a group of homeomorphisms of a locally compact σ-compact Hausdorff space X and G acts properly on X, the G\X is regular and G is equiregular and nearly proper. Below, we prove three lemmas that accomplish this task below. First we recall two definitions and establish two preliminary propositions. (These preliminary propositions are well known; we sketch their proofs for the reader's convenience.)

**Definition.** Suppose $\mathcal{M}$ is a set of maps from a topological space X to a topological space Y. If A ⊂ X and B ⊂ Y, let [A,B] = { f ∈ $\mathcal{M}$ : f(A) ⊂ B }. { [C,U] : C is a compact subset of X and U is an open subset of Y } is a subbasis for a topology on $\mathcal{M}$ called the *compact-open topology.*

**Definition.** Suppose G is a group of homeomorphisms of a topological space X. For A, B ⊂ X, let $G_{A,B}$ = { g ∈ G : g(A) ∩ B ≠ ∅ }. We say that G *acts properly on* X and that the action of G on X is *proper* if for all compact subsets A and B of X, $G_{A,B}$ is a compact subset of G when G is endowed with the compact-open topology.

**Proposition 6.1.** If G is a group of homeomorphisms of a locally compact Hausdorff space X and G acts properly on X, then G\X is Hausdorff.

**Remark.** Various forms of this proposition have appeared previously. For instance, see [Bo1] Proposition 3 in Chapter III.4.2 on page 253, [Ko] Property (ii) in Chapter 1.2 on page 3, and [P] Lemma in Section 1.2 on page 303.

**Proof.** Since locally compact Hausdorff spaces are regular ([D], Theorem XI.6.2, page 238), then any two distinct points of X have open neighborhoods with disjoint compact closures. Let a, b $\in$ X such that $\pi(a) \neq \pi(b)$. Let $A_0$ and $B_0$ be compact subsets of X such that a $\in$ int($A_0$) and b $\in$ int($B_0$). Let $\mathcal{A}$ = { A $\subset A_0$ : A is compact and a $\in$ int(A) } and let $\mathcal{B}$ = { B $\subset B_0$ : B is compact and b $\in$ int(B) }. It suffices to prove that there exists an (A,B) $\in \mathcal{A} \times \mathcal{B}$ such that $\pi$(A) $\cap \pi$(B) = $\emptyset$. Assume not. Then $G_{A,B} \neq \emptyset$ for every (A,B) $\in \mathcal{A} \times \mathcal{B}$. Also since $G_{A,B}$ = G − [A,X − B], then $G_{A,B}$ is a non-empty closed subset of $G_{A_0,B_0}$ for each (A,B) $\in \mathcal{A} \times \mathcal{B}$. Observe that if $(A_1,B_1)$, $(A_2,B_2)$, $\cdots$ , $(A_k,B_k) \in \mathcal{A} \times \mathcal{B}$, then $\cap_{1 \leq i \leq k}G_{A_i,B_i} \neq \emptyset$, because if A = $\cap_{1 \leq i \leq k}A_i$ and B = $\cap_{1 \leq i \leq k}B_i$, then (A,B) $\in \mathcal{A} \times \mathcal{B}$ and $\cap_{1 \leq i \leq k}G_{A_i,B_i} \supset G_{A,B}$. In other words, { $G_{A,B}$ : (A,B) $\in \mathcal{A} \times \mathcal{B}$ } is collection of closed subsets of the compact space $G_{A_0,B_0}$ with the *finite intersection property.* It follows that $\cap_{(A,B) \in \mathcal{A} \times \mathcal{B}}G_{A,B} \neq \emptyset$. (See [Mu], Theorem 26.9, pages 169-170.) Let g $\in$ $\cap_{(A,B) \in \mathcal{A} \times \mathcal{B}}G_{A,B}$. g(x) $\neq$ y because $\pi$(x) $\neq \pi$(y). Hence, there are open neighborhoods U of g(x) and V of Y such that cl(U) and cl(V) are compact and disjoint. Let A = $g^{-1}$(cl(U)) $\cap A_0$ and B = cl(V) $\cap B_0$. Then (A,B) $\in \mathcal{A} \times \mathcal{B}$ and g(A) $\cap$ B = $\emptyset$. Therefore, g $\notin G_{A,B}$, a contradiction. ∎

**Proposition 6.2.** Suppose G is a group of homeomorphisms of a locally compact Hausdorff space X and G is endowed with the compact-open topology. Then the action of G on X is continuous. In other words, the function $\alpha$ : G $\times$ X $\rightarrow$ X defined by $\alpha$(g,x) = g(x) is continuous.

**Proof.** Let g $\in$ G, let x $\in$ X and let U be an open neighborhood of g(x) in X. Then there is an open neighborhood V of x in X such that cl(V) is compact and g(cl(V)) $\subset$ U. Then [cl(V),U] $\times$ V is an open neighborhood of (g,x) in G $\times$ X. Moreover, if $(g',x') \in$ [cl(V),U] $\times$ V, then $\alpha(g',x')$ = $g'(x') \in$ U. Hence, $\alpha$([cl(V),U] $\times$ V) $\subset$ U. This proves $\alpha$ is continuous. ∎

**Lemma 6.3.** If G is a group of homeomorphisms of a locally compact $\sigma$-compact Hausdorff space X and G acts properly on X, then G\X is regular.

**Proof.** Proposition 6.1 implies that G\X is Hausdorff. Thus, compact subsets of G\X are closed. Since X is locally compact and $\pi$ : X $\rightarrow$ G\X is an open map, then every point of G\X lies in the interior of a compact set. Hence, every point of G\X has an open neighborhood whose closure is compact; in other words, G\X is locally compact. Since locally compact Hausdorff spaces are regular ([D], Theorem XI.6.2, page 238), then G\X is regular.

**Lemma 6.4.** If G is a group of homeomorphisms of a locally compact σ-compact Hausdorff space X and G acts properly on X, then G is equiregular.

**Proof.** Let $x \in X$ and let U be an open neighborhood of x in X. Then there is an open neighborhood V of x in X such that cl(V) is compact and cl(V) $\subset$ U. For each $y \in X$, we will find an open neighborhood $N_y$ of y in X such that for all $g \in G$, if $g(N_y) \cap cl(V) \neq \emptyset$, then $g(N_y) \subset U$. This will prove that G is equiregular.

Let $y \in X$ and let L be an open neighborhood of y in X such that cl(L) is compact. Since G acts properly on X, then $G_{cl(L),cl(V)}$ is a compact subset of G. For each $g \in G$, choose an open neighborhood $M_g$ of y in X so that cl($M_g$) is compact and:

- if $g(y) \in cl(V)$, then $g(cl(M_g)) \subset U$, and
- if $g(y) \notin cl(V)$, then $g(cl(M_g)) \subset X - cl(V)$.

Then choose $M_g{}' \in \{ U, X - cl(V) \}$ so that $g(cl(M_g)) \subset M_g{}'$. Therefore, $g \in [cl(M_g),M_g{}']$. Observe that $\{ [cl(M_g),M_g{}'] : g \in G_{cl(L),cl(V)} \}$ is a cover of the compact set $G_{cl(L),cl(V)}$ by open subsets of G. Hence, there is a finite subset F of $G_{cl(L),cl(V)}$ such that $\{ [cl(M_f),M_f{}'] : f \in F \}$ covers $G_{cl(L),cl(V)}$. Let $N_y = L \cap (\cap_{f \in F} M_f)$. Then $N_y$ is an open neighborhood of y in X. Let $g \in G$ such that $g(N_y) \cap cl(V) \neq \emptyset$. Then $g(cl(L)) \cap cl(V) \neq \emptyset$. Hence, $g \in G_{cl(L),cl(V)}$. Therefore, $g \in [cl(M_f),M_f{}']$ for some $f \in F$. Thus, $g(N_y) \subset g(M_f) \subset M_f{}'$. Since $g(N_y) \cap cl(V) \neq \emptyset$, then $M_f{}' \neq X - cl(V)$. Hence, $M_f{}' = U$. Consequently, $g(N_y) \subset U$. This proves G is equiregular. ∎

**Lemma 6.5.** If G is a group of homeomorphisms of a locally compact σ-compact Hausdorff space X and G acts properly on X, then G is nearly proper.

**Proof.** Let A and B be compact subsets of X. Since G acts properly on X, then $G_{A,B} \times A$ is a compact subset of G × X. According to Proposition 6.2, the function $\alpha : G \times X \to X$ defined by $\alpha(g,x) = g(x)$ is continuous. Hence, $\alpha(G_{A,B} \times A)$ is compact. Since $\alpha(G_{A,B} \times A) = \cup\{ g(A) : g \in G$ and $g(A) \cap B \neq \emptyset \}$, then it follows that G is nearly proper. ∎

Clearly, Theorem 1.5 and Corollary 1.6 follow from Theorem 1.3 and Corollary 1.4 using Lemmas 6.3, 6.4 and 6.5.

**References.**